\documentclass[[leqno,12pt]{amsart}
\makeatletter
\renewcommand\paragraph{\@startsection{paragraph}{4}{\z@}%
            {-2.5ex\@plus -1ex \@minus -.25ex}%
            {1.25ex \@plus .25ex}%
            {\normalfont\normalsize\bfseries}}
\makeatother
\makeatletter
\renewcommand{\@seccntformat}[1]{%
  \ifcsname format#1\endcsname
    \csname format#1\endcsname
  \else
    \csname the#1\endcsname\quad
  \fi
}
\makeatother
\usepackage{amsmath}\usepackage{MnSymbol}
\usepackage{amsfonts}\usepackage{verbatim}\usepackage{bbold}
\usepackage{graphicx}\usepackage{fontawesome}
\usepackage{marvosym}\usepackage{ifsym}
\usepackage{dictsym}\usepackage{bbm}
\usepackage{wasysym}\usepackage{dutchcal}
\usepackage{amstext}
\usepackage{amsbsy}
\usepackage{amsopn}\usepackage{hyperref}
\usepackage{amsthm}\usepackage{color}
\usepackage{pigpen}
\DeclareMathAlphabet{\mathpzc}{OT1}{pzc}{m}{it}
\def\proclaim#1{\vskip0.5em\noindent{\bf #1}\it}
\def\endproclaim{\vskip0.5em\par\noindent\rm}
\def\proclaim#1{\vskip0.5em\noindent{\bf #1}\it}
\def\endproclaim{\vskip0.5em\par\noindent\rm}
\def\demo#1{\vskip0.5em\noindent{\bf #1\ }}

\def\undersetbrace#1\to#2{\underbrace{#2}_{#1}}

\def\text#1{\mbox{#1}}
\def\flushpar{\par\noindent}

\def\tag#1{\leqno{(#1)}}
\def\mod{\mbox{ mod }}
\newcommand{\mapright}[1]{%
    \smash{\mathop{%
        \hbox to 1cm{\rightarrowfill}
        }
    \limits^{#1}
    }
}
\newcommand{\mapleft}[1]{%
    \smash{\mathop{%
        \hbox to 1cm{\rightarrowfill}
        }
    \limits_{#1}
    }
}

\def\e{\varepsilon}
\def\a{\alpha}
\def\b{\beta}
\def\g{\gamma}
\def\d{\delta}
\def\th{\theta}
\def\D{\Delta}
\def\s{\sigma}

\def\x{\times}

\def\f{\flushpar}
\def\u{\underline}
\def\v{\varphi}

\def\om{\omega}
\def\Om{\Omega}
\def\B{\mathcal B}

\def\T{\widehat T}
\def\({\biggl(}
\def\){\biggr)}
\def\<{\langle}
\def\>{\rangle}
\def\({\biggl(}
\def\){\biggr)}
\def\[{\biggl[}
\def\]{\biggr]}
\def\bul{\smallskip\f$\bullet\ \ \
$}\def\sbul{\f$\bullet\ \ \ $}\def\Par{\smallskip\f\P}
\def\pf{\smallskip\f{\tt Proof}\ \ \ \ }\def\sms{\smallskip\f}
\def\lra{\longrightarrow} \def\lfl{\lfloor}\def\rfl{\rfloor}\def\lcl{\lceil}\def\rcl{\rceil}

\let\oldtocsection=\tocsection

\let\oldtocsubsection=\tocsubsection

\let\oldtocsubsubsection=\tocsubsubsection

\renewcommand{\tocsection}[2]{\hspace{0em}\oldtocsection{#1}{#2}}
\renewcommand{\tocsubsection}[2]{\hspace{1em}\oldtocsubsection{#1}{#2}}
\renewcommand{\tocsubsubsection}[2]{\hspace{2em}\oldtocsubsubsection{#1}{#2}}
\begin{document}
\title{Extravagance, irrationality and  Diophantine approximation.}

\author{ Jon. Aaronson $\&$ Hitoshi Nakada}
\address[Aaronson]{School of Math. Sciences, Tel Aviv University
69978 Tel Aviv, Israel.}
\email{aaro@tau.ac.il}

\address[Nakada]{ Dept. Math., Keio University,Hiyoshi 3-14-1 Kohoku,
 Yokohama 223, Japan}\email[Nakada]{nakada@math.keio.ac.jp}
 \begin{abstract}  For an  invariant probability measure for the Gauss map, almost all numbers are  Diophantine if the log of the partial quotient function is integrable. We show that with respect to a  ``continued fraction mixing'' measure for the Gauss map with the log of the partial quotient function  non-integrable, almost all numbers are Liouville.  We
  also exhibit Gauss-invariant, ergodic measures with arbitrary  irrationality exponent.
  \

  The proofs are applications of our study of the  ``extravagance'' of positive, stationary, stochastic processes.
  \

In addition, we prove a Khinchin-type dichotomy for Diophantine approximation with respect to  ``weak Renyi measures'' which are ``doubling at $0$''.
 \end{abstract}
\thanks{\copyright 2023-24.}\dedicatory{ Dedicated to the memory of Yuji Ito.}
\maketitle\markboth{\copyright J. Aaronson and H. Nakada}{Extravagance, irrationality and  Diophantine approximation.}
\tableofcontents
\setcounter{tocdepth}{2}
\section*{\S1 Introduction}

\subsection*{Stationary processes}
\

A {\it stochastic process}  with values in a measurable space $Z$ is a quadruple
$(\Om,m,\tau,\Phi)$ where $(\Om,m,\tau)$ is a non-singular transformation  and  $\Phi:\Om\to Z$ is measurable.

\

It is
\sbul\ {\it forward generating} if
$\s(\{\Phi\circ\tau^k:\ k\ge 0\})\overset{m}=\B(\Om)$;
\

\sbul {\it stationary} if $(\Om,m,\tau)$ is a probability preserving  transformation and
\

\sbul\ {\it ergodic}  if $(\Om,m,\tau)$ is an ergodic probability preserving transformation.
\

\subsection*{Partial quotients}
\

Let
$\mu\in\mathcal{P}(\Bbb I)$  be invariant  under the  {\it Gauss map}  $G:\mathbb{I}:=[0,1]\setminus\Bbb Q\hookleftarrow$, defined by
$$G(x):=\{\tfrac1x\}=\tfrac1x-\lfl\tfrac1x\rfl.$$

As shown in in \cite{Kh}, various
  Diophantine properties of $\mu$-typical $x\in\Bbb I$ are determined
by the asymptotic properties of the (stationary) process of {\tt partial quotients} $(\Bbb I,\mu,G,a)$ (with
 $a(x):=\lfl \tfrac1x\rfl$).
 \

 In this situation, we'll  consider
  the {\tt extravagance} of the stationary process $(\Bbb I,\mu,G,\log a)$.

\

 \subsection*{Extravagance}\ \

The {\it extravagance}  of the non-negative sequence $(x_n:\ n\ge 1)\in[0,\infty)^{\Bbb N}$ is
$$\mathbbm{e}((x_n:\ n\ge 0)):=\varlimsup_{n\to\infty}\frac{x_{n+1}}{\sum_{k=1}^{n}x_k}\in [0,\infty]$$ if $\exists\ n\ge 1,\ x_n>0;\ \&\ \mathbbm{e}(\overline{0}):=0$.

The {\it extravagance}  of the non-negative stationary process $(\Om,m,\tau,\Phi)$ is the random variable
$\mathfrak{e}(\Phi,\tau)$ on $(\Om,m)$ defined by
$$\mathfrak{e}(\Phi,\tau)(\om):=\mathbbm{e}((\Phi(\tau^n\om):\ n\ge 0)).$$

\

Calculation shows that
$ \mathfrak{e}(\Phi,\tau)\circ\tau\ge \mathfrak{e}(\Phi,\tau)$ a.s. and the extravagance is a.s. constant if $(\Om,m,\tau)$ is ergodic.
\

It follows from the ergodic theorem that for a stationary process, $\Bbb E(\Phi)<\infty\ \Rightarrow\ \mathfrak{e}(\Phi,\tau)=0$ a.s..

\

We show (Theorem 2.2 on p.\pageref{2.2}) that if the non-negative stationary process $(\Om,m,\tau,\Phi)$ is {\tt  continued fraction mixing} (i.e. satisfies \ref{CF} on p.\pageref{CF}), then $\mathfrak{e}(\Phi,\tau)=0$
a.s. iff $\Bbb E(\Phi)<\infty$ and otherwise
$\mathfrak{e}(\Phi,\tau)=\infty$ a.s..
\

On the other hand,
\sbul there is a {\tt Markov shift} $(\Bbb N^\Bbb Z,m,S=\text{\tt shift})$ so that
for any $r\in \Bbb R_+$
 there is a
{\tt finitary} function $\Phi^{(r)}:\Bbb N^\Bbb Z\to\Bbb R_+$ with $\mathfrak{e}(\Phi^{(r)},S)=r$ a.s. (Theorem 3.1 on p.\pageref{3.1});
\

and

\sbul for any aperiodic, ergodic, probability preserving transformation
$(X,m,T)$, for any $r\in \Bbb R_+$
there is a measurable function $g^{(r)}:X\to\Bbb R_+$ so that $\mathfrak{e}(g^{(r)},T)=r$ a.s. (Theorem 3.2 on p.\pageref{3.2}).

 \subsection*{Irrationality}\label{irr}
 Let $\mathbb{I}:=[0,1]\setminus\Bbb Q$ be the irrationals in $(0,1)$.
\

An irrational $x\in\mathbb{I}$ is called {\it badly approximable of order} $s>0$ (abbr. $s$-{\tt BA}) if
 $\min_{0\le p\le q}|x-\tfrac{p}q|\gg\tfrac{1}{q^s}$ as $q\to\infty$.
 \

  The {\it irrationality} (exponent) of $x\in\mathbb{I}$  (as in \cite[Appendix E]{Bugeaud}) is
$$\mathbbm{i}(x):=\inf\,\{s>0:\ x\ \text{is}\  s-\text{\tt BA}\}\le\infty.$$
\

 By Dirichlet's theorem, $\mathbbm{i}\ge 2$.

\

By Legendre's theorem (see e.g. \cite[Theorem 5C]{WSchmidt}),
for $x\in\Bbb I$, if $p,q\in\Bbb N,\ \text{\tt\small gcd}(p,q)=1$ and
$|\tfrac{p}q-x|<\tfrac1{2q^2}$, then $\tfrac{p}q=\tfrac{p_n(x)}{q_n(x)}$ (some $n\ge 1$) where $(\tfrac{p_n(x)}{q_n(x)}:\ n\ge 1)$ are the {\tt convergents} of $x$ (as on p.\pageref{cgts}).

\

It follows that $x\in\mathbb{I}$ is $s$-{\tt BA} ($s\ge 2$) iff $|x-\tfrac{p_n(x)}{q_n(x)}|\gg\tfrac{1}{q_n(x)^s}$ as $n\to\infty$, whence
\begin{align*}\tag*{\faBullseye}\label{faBullseye}
\mathbbm{i}(x)=\inf\,\{s>2: \ |x-\tfrac{p_n(x)}{q_n(x)}|\gg\tfrac{1}{q_n(x)^s}\ \text{as}\ n\to\infty\}.
\end{align*}

An irrational $x\in\mathbb{I}$  is called
\sbul {\it Diophantine}\ if $\mathbbm{i}(x)=2$;
\sbul {\it very well approximable} if $\mathbbm{i}(x)>2$; and
\sbul  {\it Liouville} if $\mathbbm{i}(x)=\infty$.
\

It is shown in \cite{BugExact} that for $s\ge 2$, the Hausdorff dimension of the set $\{x\in\Bbb I:\ \mathbbm{i}(x)=s\}$ is $\tfrac{2}s$.

\

It turns out that (Bugeaud's Lemma  on page \pageref{4.4}) for $x\in\Bbb I$,
\begin{align*}\tag*{\ref{faKey}}
\mathbbm{i}(x)=2+\mathbbm{e}&((\log \tfrac1{G^n(x)}:\ n\ge 0)).
\end{align*}
and for $G$-invariant $\mu\in\mathcal{P}(\Bbb I)$:
\begin{align*}\tag*{\ref{faPaw}}
\mathbbm{i}=2+\mathfrak{e}(\log a,G)\ \
\mu-\text{a.s.;}
\end{align*}
whence if
$\Bbb E_\mu(\log a)<\infty$, then $\mu$-a.s., $\mathfrak{e}(\log a,G)=0$ and
$$\mathbbm{i}=2+\mathfrak{e}(\log a,G)=2.$$
By Corollary 4.5 (on p.\pageref{4.5}): if $\mu\in\mathcal{P}(\Bbb I)$ is so that
$(\Bbb I,\mu,G,a)$ is stationary and  continued fraction mixing, then
\sbul if $\Bbb E_\mu(\log a)<\infty$, then $\mu$-a.e. $x\in\Bbb I$ is Diophantine; and
\sbul if $\Bbb E_\mu(\log a)=\infty$, then $\mu$-a.e. $x\in\Bbb I$ is Liouville;
\sbul $\forall\ r\ge 2,\ \exists\ \mu\in\mathcal{P}(\Bbb I)$ so that $(\Bbb I,\mu,G,a)$ is an ergodic, stationary process and so that $\mathbbm{i}=r\ \mu$-a.s..

\subsection*{A Khinchin-type dichotomy for $G$-invariant measures}
\

It is shown in \cite{RenRep,Fexp} that {\it Gauss measure} $\mu\in\mathcal{P}(\Bbb I),\ d\mu(x)=\tfrac{dx}{\log 2(1+x)}$ is a  Renyi measure for $G$ in that $(\Bbb I,\mu,G,a)$
has the {\tt Renyi property} (as in \ref{Rp} on p.\pageref{Rp}) and in \cite{GM} it is shown  that $(\Bbb I,\mu,G,a)$ is a Gibbs-Markov map whence  {\tt continued fraction mixing} (as in \ref{CF} on p.\pageref{CF}).
\

We'll call a $G$-invariant measure $\nu\in\mathcal{P}(\mathbb{I})$:

{\it Renyi, weak Renyi} or {\it continued fraction mixing}  according to whether the stationary process $(\Bbb I,\nu,G,a)$ has the
{\tt Renyi property} (as in \ref{Rp}), the {\tt weak Renyi property} (as in \ref{wRp}), or is {\tt continued fraction mixing} (as in \ref{CF}); respectively.

In \S5 we establish a Khinchin type dichotomy for ergodic, weak Renyi measures  which are    {\tt doubling at $0$} as in \ref{faCopy}
 (Theorem 5.1 on p.\pageref{5.1}).

\subsection*{Renyi properties and continued fraction mixing}
\

The  stationary, forward generating, stochastic process  $(\Om,m,\tau,\Phi)$
\bul is {\it independent} if $\{\Phi\circ\tau^n:\ n\ge 1\}$ are independent random variables;
\bul has the  {\it  Renyi property} if
\begin{align*}\tag{$\overline{\mathfrak{R}}$}\label{Rp}\exists\ & M>1\ \text{s.t.}\
m(A\cap B)=M^{\pm 1}m(A)m(B)\ \forall\ n\ge 1,\ \\ & A\in\s(\{\Phi\circ\tau^k:\ 0\le k\le n\}),\ B\in
 \s(\{\Phi\circ\tau^\ell:\ \ell\ge n+1\});
\end{align*}
\bul has the  {\it  weak Renyi property} if
\begin{align*}\tag{$\underline{\mathfrak{R}}$}\label{wRp}\exists\ & M>1\ \text{s.t.}\
m(A\cap B)\le Mm(A)m(B)\ \forall\ n\ge 1,\ \\ & A\in\s(\{\Phi\circ\tau^k:\ 0\le k\le n\}),\ B\in
 \s(\{\Phi\circ\tau^\ell:\ \ell\ge n+1\});
\end{align*}
\bul is {\it continued fraction {\rm (abbr. {\tt  c.f.})} mixing} if $\exists\
(\vartheta(N):\ N\ge 1)\in\Bbb R_+^\Bbb N,\  \vartheta(N)\downarrow 0$ so that
\begin{align*}\tag{{\tt CF}}\label{CF}
|m(A&\cap B)-m(A)m(B)|\le \vartheta(N)m(A)m(B)\ \forall\ n\ge 1,\ \\ & A\in\s(\{\Phi\circ\tau^k:\ 0\le k\le n\}),\ B\in
 \s(\{\Phi\circ\tau^\ell:\ \ell\ge n+N\}).
\end{align*}

\

Note that a {\tt c.f.} mixing process has the weak Renyi property, but not necessarily the Renyi property.
For example, a stationary, mixing  {\tt Gibbs-Markov map} $(X,m,T,\a)$  (as in \cite{GM})  is  weak Renyi,  but has the   Renyi property  if and only if $Ta=X\ \forall\ a\in\a$.

As shown in \cite{RenRep}:
 a stationary,  Renyi process $(X,m,T,\Phi)$  is {\it exact}
in the sense that the {\it tail} field is trivial:
$$\mathcal{T}(T):=\bigcap_{n\ge 1}T^{-n}\B(X)\overset{m}=\{\emptyset,X\}.$$
It follows from \cite[Theorem 1]{Bra} that a stationary process with the
Renyi property  is {\tt c.f.} mixing.

\

 A stationary, weak Renyi process $(X,m,T,\Phi)$  need not be ergodic. For example if $(X,m,T,\Phi)$ is an $\Bbb N$-valued Renyi process, then $(X\x\{0,1\},m\x\#,T\x\text{\tt Id},\tilde{\Phi})$ (with $\tilde{\Phi}(x,y):=\Phi(x)+\sqrt 2y$) is weak Renyi but not ergodic.
 \

 However, a stationary, weak Renyi process $(X,m,T,\Phi)$  has a finite  tail field and  hence is  exact if totally ergodic.
\

To see that $\mathcal{T}(T)$ is purely atomic, let $A\in\mathcal{T}(T),\ m(A)>0$
and let $A_n\in\s(\{\Phi\circ T^k:\ 0\le k<n\})$,\ $m(A_n\D A)\xrightarrow[n\to\infty]{}0$, then,
\begin{align*}
m(A)&\xleftarrow[n\to\infty]{}m(A_n\cap A)=m(A_n\cap T^{-n}T^nA)
\le Mm(A_n)m(T^nA)\ \text{ by \ref{wRp}}\\ &=Mm(A_n)m(A)\ \xrightarrow[n\to\infty]{}Mm(A)^2
\end{align*}

and $m(A)\ge \tfrac1M$. Thus  $\#\mathcal{T}(T)<\infty$ and  the Pinsker (i.e. tail)  factor  consists of finitely many periodic, ergodic components.
Thus, $T$ is exact 
 if  totally ergodic.

\subsection*{Fibered systems}
\

As in \cite{Schweiger}, a (stationary) {\it fibered system}  $(X,m,T,\a)$  is  a probability preserving transformation $T$  of a standard probability
space  $(X,m)$, equipped with a countable (or finite), measurable partition $\a$
which generates $\B(X)$ under $T$ in the sense
that $\s(\{T^{-n}\a:\ n\ge 0\})=\B$ and which  satisfies
 $T:a\to Ta$ invertible
and nonsingular for $a\in\a$.
\

A fibered system $(X,m,T,\a)$ can also be viewed as a forward generating, stochastic process $(X,m,T,\Phi)$ with
$\Phi:X\to\a,\ x\in\Phi(x)\in\a$ and we call it  {\it Renyi, weak Renyi} or {\tt c.f.}{\it mixing}
accordingly.


\section*{\S2 Extravagance  of continued fraction mixing  processes}

\

\proclaim{2.1 Proposition}\label{2.1}
\

\ \  Let $(\Om,m,\tau,\Phi)$ be a  stationary process. Suppose that $f:\Om\to [0,\infty),\ \Bbb E(f)<\infty$, then $m$-a.s.:
$$\mathfrak{e}(\Phi+f,\tau)=\mathfrak{e}(\Phi,\tau).$$

\endproclaim\demo{Proof}\ \ There is no loss in generality in assuming that $\tau$ is ergodic and that $\Bbb E(f),\ \Bbb E(\Phi)> 0$.
\

If $\Bbb E(\Phi)<\infty$, then $\Bbb E(\Phi+f)<\infty$ and
$$\mathfrak{e}(\Phi+f,\tau)=\mathfrak{e}(\Phi,\tau)=0.$$
\

Now suppose that $\Bbb E(\Phi)=\infty$.

By the ergodic theorem, writing $g_n^{(\tau)}:=\sum_{k=0}^{n-1}g\circ\tau$ for $g=f,\Phi$,
$$\tfrac{f_n^{(\tau)}}n\xrightarrow[n\to\infty]{}\Bbb E(f),\ \tfrac{\Phi_n}n\xrightarrow[n\to\infty]{}\infty\ \text{$m$-a.e. .}$$ Moreover $f\circ\tau^n=o(n)$ a.s.,
whence
$\tfrac{(f+\Phi)\circ\tau^n}{(f+\Phi)^{(\tau)}_n} \sim\tfrac{\Phi\circ\tau^n}{\Phi^{(\tau)}_n}$ and $$\mathfrak{e}(\Phi+f,\tau)=\mathfrak{e}(\Phi,\tau).\ \ \CheckedBox$$
\

\proclaim{2.2 Theorem}\label{2.2}\ \
\

Suppose that $(\Om,\mu,\tau,\a)$ is a continued fraction mixing, probability preserving fibered system and that
$\Phi:\Om\to\Bbb N$ is $\a$-measurable, then

$$\mathfrak{e}(\Phi,\tau)=\begin{cases}& 0 \ \ \ \ \text{\rm a.s.\ \ if} \ \ \ \ \ \ \Bbb E(\Phi)<\infty\ \ \ \&\\ &
 \infty \ \ \  \text{\rm a.s.\ \ if} \ \ \ \ \  \Bbb E(\Phi)=\infty.
               \end{cases} $$
\endproclaim

 In the independent case the result is proved in \cite{Raugi} (see also \cite{C-Z} for related results).

\

 The proof of Theorem 2.2 involves
 \subsection*{Kakutani skyscrapers}\label{skyscraper}
 \

 Let $(\Om,\mu,\tau,\phi)$ be a $\Bbb N$-stationary process.
 \

 The   {\it Kakutani skyscraper}
(as in   \cite{Kak})  is the conservative, ergodic, measure preserving transformation  $(\Om,\mu,\tau)^\phi:=(X,m,T)$ where
\begin{align*}\tag*{\faBuildingO}\label{faBuildingO}
 X:=\{(\om,n)\in\Om &\x\Bbb N:\ 0\le n\le \phi(\om)-1\},\ m:=(\mu\x\#)|_X\ \&\\ & T(\om,n):=\begin{cases}& (\om,n+1)\ \ n<\phi(\om)-1\\ &   (\tau(\om),1)\ \ n=\phi(\om)-1.                                                                                \end{cases}
\end{align*}
\subsection*{Renewal Process}\label{renewal}
\

A {\it renewal process} is a Kakutani skyscraper $(\Om,\mu,\tau)^\phi$ where $(\Om,\mu,\tau,\phi)$  is independent. It is isomorphic to the Markov shift with state space $\Bbb N$ and transition matrix given by
$$p_{s,t}=\begin{cases}& \mu([\phi=t])\ \ \ s=1;\\ & 1\ \ \ s=t+1;\\ & 0\ \ \ \text{else;}
          \end{cases}$$
with stationary distribution $\rho\in\mathfrak{M}(\Bbb N)$ given by $\rho_t=\mu([\phi\ge t])$.
\

That is $(X,m,T)=(\Om,\mu,\tau)^\phi\cong (\Bbb N^\Bbb Z,m_{p,\rho},\text{\tt shift})$
where
$$m_{p,\rho}([s_0,\dots,s_n]_k)=\rho_{s_0}p_{s_0,s_1}\dots p_{s_{n-1},s_n}$$
with $[s_0,\dots,s_n]_k:=\{x\in \Bbb N^\Bbb Z:\ x_{k+j}=s_j\ \forall\ 0\le j\le n\}$.
\

By ergodicity and recurrence of $(\Om,m,S)$, for a.e. $\om\in\Om$,
$$K(\om):=\{n\in\Bbb Z:\ \om_n=1\}$$ has no infinite gaps ($[\a,\b]\subset \Bbb Z\setminus K(\om)\ \Rightarrow\ \b-\a<\infty$).

Write $K(\om):=\{c_n(\om):\ n\in\Bbb Z\}$ where $c_0\le 0<c_1$.

The isomorphism $\mathcal{c}:(\Bbb N^\Bbb Z,m_{p,\rho},\text{\tt shift})\to(X,m,T)=(\Om,\mu,\tau)^\phi$ is given by the correspondence
\begin{align*}\tag*{\faArrowsH}\label{faArrowsH}\begin{split}&
 \om\in\Bbb N^\Bbb Z\leftrightarrow \mathcal{c}(\om)=(\eta,\mathcal{k})\in X\subset \Bbb N^\Bbb Z\x\Bbb N
  \ \ \text{ where}\\ &\eta(\om)=(\eta_n(\om)=c_{n+1}(\om)-c_n(\om):\ n\in\Bbb Z)\ \&\ \mathcal{k}(\om)=-c_0(\om)+1.
                                                \end{split}
\end{align*}

 \subsection*{Darling-Kac sets}
 \

A {\it Darling-Kac set} (as in \cite{Da-Kac}) for the measure preserving transformation $(X,m,T)$ is a set $A\in\B(X),\ 0<m(A)<\infty$ so that
$$\frac1{a_n(A)}\sum_{k=0}^{n-1}\T^k1_A\xrightarrow[n\to\infty]{}\ m(A)$$
uniformly on  $A$\ with\ $a_n(A):=\sum_{k=0}^{n-1}\tfrac{m(A\cap T^{-k}A)}{m(A)^2}.$
\

If the conservative, ergodic, measure preserving transformation $(X,m,T)$ has a Darling-Kac set $A$, then $T$ is {\it pointwise dual ergodic} in the sense that

there is a sequence $a(n)=a_n(T)$ (the {\it return sequence} of $(X,m,T)$) so that
\begin{align*}\tag{{\tt PDE}}\label{PDE}\frac1{a(n)}\sum_{k=0}^{n-1}\T^kf\xrightarrow[n\to\infty]{}\ \int_Xfdm\ \text{a.e.}\ \forall\ f\in L^1(m).
\end{align*}
Here $\T:L^1(m)\hookleftarrow$ is the {\it transfer operator}   defined  by
        $$\int_A\T fdm\ =\ \int_{T^{-1}A} fdm\ \ \ A\in\B(X)$$ and
        $a_n(A)\sim a_n(T)$ for any Darling-Kac set $A$.
        See \cite{Aasdist} (also \cite[\S3.7]{A1})
        \


Let $(\Om,m,\tau,\a)$ be an ergodic,  probability preserving fibered system and let
$\Phi:\Om\to\Bbb N$ be $\a$-measurable.
We'll need the following facts about the Kakutani skyscraper $(X,m,T)=(\Om,m,\tau)^\Phi$.
\Par1 If $(\Om,m,\tau,\a)$ is continued fraction mixing, then  $\Om$ is a Darling-Kac set for $T$. See \cite{Arfx} (and \cite{Da-Kac} for the independent case).
\Par2 If $\Om$ is a Darling-Kac set for $T$, then
\begin{align*}\tag*{\dsjuridical}\label{dsjuridical}
a_n(T)=2^{\pm 1}\overline{a}(n)\ \text{where}\ \overline{a}(n):=\tfrac{n}{L(n)}\ \text{with}\ L(n):=\Bbb E(\Phi\wedge n).
\end{align*}
See \cite[Theorem 3]{Aasdist} (also \cite[Lemma 3.8.5]{A1}). Note that \ref{dsjuridical} is an elementary consequence of the {\tt discrete renewal equation} as in
\cite[\S1.8]{Chu} in the independent case.
We'll need
\proclaim{Lemma 2.3}\label{2.3}\ \ Let $\xi$ be an $\Bbb N$-valued random variable.
\

If $\Bbb E(\tfrac{\xi}{L(\xi)})<\infty$ with $L(t):=\Bbb E(\xi\wedge t)$, then $\Bbb E(\xi)<\infty$.
\endproclaim

\demo{Proof}\ Let $(\Om,\mu,\s):=(\Bbb N^\Bbb N,\text{\tt\small dist}(\xi)^\Bbb N,\text{\tt\small) shift})$ and define $\Phi:\Om\to\Bbb N$ by $\Phi(\om):=\om_1$, then
with $\xi_n:=\Phi\circ\s^{n-1}$, $\xi_n:\ n\ge 1)$ are independent, identically distributed random variables each distributed as $\xi$.

Let $(X,m,T):=(\Om,\mu,\s)^\Phi$.
\

By \P1 $\Om$ is a Darling-Kac set for $T$ and by \P2,  $a_n(T)=2^{\pm 1}\overline{a}(n)$ with $\overline{a}(n):=\tfrac{n}{L(n)}$.
\

Now suppose that the lemma fails and  $\Bbb E(\overline{a}(\xi))<\infty$  whereas
$\Bbb E(\xi)=\infty$.

Now $m(X)=\Bbb E(\xi)=\infty$ entails $\tfrac{\overline{a}(n)}n\ \downarrow 0$ whence by
Feller's theorem (\cite{Fe1})
\begin{align*}
 \tag*{\faBug}\label{faBug}\tfrac{\Phi_n}{b(n)}\xrightarrow[n\to\infty]{}0\ \text{a.s.}\
\end{align*}
with $\Phi_n:=\sum_{k=0}^{n-1}\Phi\circ\s^k=\sum_{k=1}^{n}\xi_k$ and $b:=\overline{a}^{-1}$.

It follows from this (\cite{A8} -- also \cite[Theorem 2.4.1]{A1}) that
$$\tfrac1{\overline{a}(n)}\sum_{k=0}^{n-1}1_\Om\circ T^k\xrightarrow[n\to\infty]{}\infty\ \text{a.s.}$$
whence by Fatou's lemma
$$2\ge \tfrac{a_n(T))}{\overline{a}(n)}=\int_\Om(\tfrac1{\overline{a}(n)}\sum_{k=0}^{n-1}1_\Om\circ T^k)dm\xrightarrow[n\to\infty]{}\infty.
\ \ \XBox$$
Thus $\Bbb E(\Phi)<\infty$.\ \ \CheckedBox
\demo{Proof of Theorem 2.2}
\

As mentioned above,  $\Bbb E(\Phi)<\infty\ \Rightarrow\ \mathfrak{e}(\Phi,\tau)=0$ a.s. by the ergodic theorem. It suffices to prove that $\mathfrak{e}(\Phi,\tau)<\infty\ \Rightarrow\ \Bbb E(\Phi)<\infty$ for which, by  Lemma 2.3, $\Bbb E(a(\Phi))<\infty$ suffices.

\

Assume $\mathfrak{e}(\Phi,\tau)<\infty$ a.s..

We show first that $\exists\ \g\in\Bbb N$ so that
\begin{align*}\tag*{\Wheelchair}\label{Wheelchair}\sum_{n\ge 1}\mu([\Phi\circ \tau^n>\g\Phi_n])<\infty.
\end{align*}
\pf of \ \ref{Wheelchair}
\

For $\d>0$ set $A_n(\d):=[\Phi\circ \tau^n>\d\Phi_n]\in\s(\a_{n+1})$,
 then for $n,\ k\ge 2$
  \begin{align*}
   A_n(\d)\cap A_{n+k}(\d)&=
[\Phi\circ\tau^n>\d\Phi_n\ \&\ \Phi\circ\tau^{n+k}>\d\Phi_{n+k}]\\ &
\subseteq [\Phi\circ\tau^n>\d\Phi_n\ \&\ \Phi\circ\tau^{n+k}>\d\Phi_{k-1}\circ\tau^{n+1}]\\ &= A_n(\d)\cap\tau^{-(n+1)}A_{k-1}(\d)
  \end{align*}

whence by the weak Renyi property (entailed by continued fraction mixing),
$$\mu(A_n(\d)\cap A_{n+k}(\d))\le M\mu(A_n(\d))\mu(A_{k-1}(\d)).$$
Thus, with $N_n:=\sum_{k=1}^n1_{A_k(\d)}$,
  \begin{align*}\tag*{\faTrophy}\label{faTrophy}\Bbb E(N_n^2)\le 3\Bbb E(N_n)+2M\Bbb E(N_n)^2.
  \end{align*}

Fix $\eta>\mathfrak{e}(\Phi,\tau)$, then $\sum_{n\ge 1}1_{A_n(\eta)}<\infty$ a.s.
\
By \ \ref{faTrophy}\ \
and the Erdos-Renyi Borel-Cantelli lemma (\cite{ErdosRenyi}\ $\&$/or  \cite[p.391]{Ren2})
$$\sum_{n\ge 1}\mu(A_n(\eta))<\infty.\ \ \CheckedBox\ \text{\ref{Wheelchair}}$$
\

Let $(X,m,T)=(\Om,\mu,\tau)^\Phi$ be the  Kakutani skyscraper  as in \ref{faBuildingO}.
\

By \P1 (p.\pageref{dsjuridical}), $(X,m,T)$ is a
pointwise dual ergodic measure preserving transformation with
$$a_n(T)=a(n)=\sum_{k=0}^{n-1}m(\Om\x\{1\}\cap T^{-k}\Om\x\{1\})$$
and $\Om\x\{1\}$ is a  Darling-Kac set for $T$.
\

Thus, by \P2 (p.\pageref{dsjuridical}), $\exists\ M>1\ \&\ N_0\in\Bbb N$ so that

\begin{align*}\tag*{\dstechnical}\label{dstechnical}
s_n:=\sum_{k=1}^n\T^k1_\Om\x\{1\}= &M^{\pm 1}\overline{a}(n)\ \text{on}\ \Om\x\{1\}\ \forall\ n\ge N_0
\end{align*}
where $\overline{a}(n)=\tfrac{n}{\Bbb E(\Phi\wedge n)}$ is as in \ref{dsjuridical} (p.\pageref{dsjuridical}).
\

 Finally, we claim  that
\begin{align*}\tag*{\scriptsize\faShip}\label{faShip}\Bbb E(\overline{a}(\Phi))<\infty.
\end{align*}

\pf Let $\g\in\Bbb N$ be as in \ref{Wheelchair} (p.\pageref{Wheelchair}), then
\begin{align*}\tag*{\faLeaf}\label{faLeaf}\infty&>C:=\sum_{n\ge 0}\mu([\Phi\circ \tau^n>\g\Phi_n])=
\sum_{k\ge n\ge 1}\mu([\Phi_n=k]\cap  \tau^{-n}[\Phi\ge \g k])\\ &=
\sum_{k=1}^\infty m(\Om\x\{1\}\cap T^{-k}([\Phi\ge \g k])
=\int_{\Om\x\{1\}}\sum_{k=1}^{\lfl\frac{\Phi}\g\rfl}
1_{[\Phi\ge \g k]}\T^k1_{\Om\x\{1\}} dm\\ &\ge
\int_{[\Phi\ge \g N_0]}\sum_{k=1}^{\lfl\frac{\Phi}\g\rfl}
\T^k1_{\Om\x\{1\}} dm\ge\tfrac1M\Bbb E(1_{[\Phi\ge \g N_0]}\overline{a}(\tfrac{\Phi}\g))\ \text{by\ \ref{dstechnical} on p.\pageref{dstechnical}.}
 \end{align*}

Using \ref{faLeaf},
 \begin{align*}\Bbb E(\overline{a}(\Phi))&\le \g\Bbb E(\overline{a}(\tfrac{\Phi}\g))\le \overline{a}(\tfrac{N_0}\g)+\g\Bbb E(\overline{a}(\tfrac{\Phi}\g) 1_{[\Phi\ge \g N_0]})\\ &\le \overline{a}(\tfrac{N_0}\g)+
 M\g\int_{\Om\x\{1\}}\sum_{k\ge 1}1_{[\Phi\ge \g k]\x\{1\}}\T^k1_{\Om\x\{1\}} dm\\ &\le \overline{a}(\tfrac{N_0}\g)+M\g C<\infty.\ \ \CheckedBox\ \text{\ref{faShip}}
 \end{align*}
This proves Theorem 2.2. \ \CheckedBox

\

\section*{\S3 Extravagance of ergodic, stationary processes}\label{s3}

Next, we obtain ergodic stationary processes with arbitrary extravagance.

\

\proclaim{3.1 Theorem}\label{3.1}
\

There is a  Markov shift\  $(\Om=\Bbb N^\Bbb Z,m,S=\text{\tt shift})$\
so that for each $t\in\Bbb R_+$ there is a finitary function $g=g_t:\Bbb N^\Bbb N\to\Bbb R_+$ so that
$\mathfrak{e}(g,S)=t$ a.s.\endproclaim
Here, a measurable function $f:\Bbb N^\Bbb Z\to\Bbb R$ is  {\it finitary} if $\exists\ N:\Bbb N^\Bbb Z\to\Bbb N\cup\{\infty\}$ measurable so that for $m$ a.e. $\om\in\Om$,
$$N(\om)<\infty\ \&\ f([\om_{-N(\om)},\dots,\om_{N(\om)}]_{-N(\om)})=\{f(\om)\}.$$
Here, for $j,k,L\in\Bbb Z,\ j<k$,
$$[a_j,a_{j+1},\dots,a_k]_L:=\{x\in \Bbb N^\Bbb Z:\ x_{L+i}=a_{j+i}\ \forall\ 0\le i\le k-j\}.$$
 \proclaim{3.2 Theorem}\label{3.2}

\ Let $(X,m,T)$ be an aperiodic, ergodic, probability preserving transformation.
\

For each $r\in \Bbb R_+,\ \exists$ an $\Bbb R_+$-valued measurable function $g=g_r:\Om\to\Bbb R_+$ so that
$$\mathfrak{e}(g,T)=r\ \text{a.s.}$$
\endproclaim

\proclaim{3.3 Main Lemma}\label{3.3}\ \ Suppose that $a>1\ \&\ (Y,p,\s,\phi)$ is an $\Bbb N$-valued, ergodic stationary process so that
\begin{align*}
 &\tag{i}\ \Bbb E(\phi)<\infty;\\ &\tag{ii}\label{ii}\ \mathfrak{e}(\sqrt a^\phi,\s)=\infty\ \text{a.s.}.
\end{align*}
\

Let $(\Om,\mu,\tau):=
(Y,\tfrac1{\Bbb E(\phi)}\cdot p,\s)^\phi$ and define $\Psi:\Om\to\Bbb R_+$ by
$$\Psi(y,n):=a^{n\wedge(\phi(y)-n)},\ \ \ (y,n)\in\Om=\{(x,\nu):\ x\in Y,\ 0\le \nu<\phi(x)\},$$
then $\mathfrak{e}(\Psi,\tau)=a-1\ \text{a.s.}.$\endproclaim\demo{Proof}\  \ For $y\in Y$, let
$$B(y):=((\Psi(\tau^m(y,0)):\ 0\le n<\phi(y)),$$
then
\begin{align*}
B(y)=
(1,a,a^2,\dots,a^{\lfl\phi(y)/2\rfl},a^{\lfl\phi(y)/2\rfl-1},\dots,a)
\end{align*}
whence $\Psi\circ\tau=a^{\pm 1}\Psi$ and
\begin{align*}\tag*{{\scriptsize\faAnchor}}\label{faAnchor}
\widetilde{\Psi}(y):=\sum_{j=0}^{\phi(y)-1}\Psi(\tau^j(y,0))=\tfrac{a+1}{a-1}\cdot(a^{\lfl\phi(y)/2\rfl}-1).
\end{align*}

\

Moreover, for fixed $y\in Y$,

 \begin{align*}\Psi^{(\tau)}_{\phi_K}(y,0)=\widetilde{\Psi}_K^{(\s)}(y).
\end{align*}
Next, for a.e. $y\in Y$,  each $n\ge 0$ has the decomposition
\begin{align*}\tag*{{\scriptsize\faBell}}\label{faBell}
&n=\phi^{(\tau)}_{K_n(y)}(y)+r_n(y)\  \text{where}\\ &
K_n(y):=\sum_{j=1}^n1_Y\circ\tau(y,0)=\#\,\{k\ge 1:\ \phi_k\le n\}\\ &\ \  \ \&\ 0\le r_n(y)<\phi(\s^{K_n}(y)).\end{align*}
Consequently,
 \begin{align*}\Psi^{(\tau)}_n(y,0)&=\Psi^{(\tau)}_{\phi_{K_n}}(y,0)+\Psi^{(\tau)}_{r_n}(\s^{K_n}y,0)\\ &=
 \widetilde{\Psi}_{K_n}^{(\s)}(y)+\Psi^{(\tau)}_{r_n}(\s^{K_n}(y,0).
 \end{align*}
Thus
 \begin{align*}\tag*{{\scriptsize\faBank}}\label{faBank} M_n(\Psi,\tau)(y,0)=\frac{\Psi(\tau^n(y,0))}{\Psi^{(\tau)}_n(y,0)}=
 \frac{a^{r_n\wedge(\Psi(\s^{K_n}y)-r_n)}}{ \widetilde{\Psi}_{K_n}^{(\s)}(y)+\Psi^{(\tau)}_{r_n}(\s^{K_n}y,0)}.
   \end{align*}
By ergodicity, it suffices to show that $\overline{M}:=\varlimsup_{n\to\infty}M_n= a-1$ a.s. on $Y$.
\demo{ Proof that  $\overline{M}\ge a-1$}
\

By \ref{ii} and \ \ref{faAnchor},\ \
$\mathfrak{e}(\widetilde{\Psi},\s)=\infty$ a.s. on $Y$.
\

For any $\e>0,\ J\ge 1$ $\&\ y\in Y$ s.t. $\mathfrak{e}(\widetilde{\Psi},\s)(y)=\infty$, $\exists\ N>J$ so that
$$a^{\lfl\phi(\s^Ny)/2\rfl}>\tfrac1\e\widetilde{\Psi}^{\s)}_N(y).$$
Let $n:=\phi_N(y)+\lfl\phi(\s^Ny)/2\rfl$, then
 \begin{align*}M_n(\Psi,\tau)(y,0)&=\frac{a^{\lfl\phi(\s^Ny)/2\rfl}}{\widetilde{\Psi}_{N}^{(\s)}(y)+
  \Psi^{(\tau)}_{\lfl\phi(\s^Ny)/2\rfl}(\s^{N}y,0)}\ \ \text{by \ \ref{faBank}}\\ &=
  \frac{a^{\lfl\phi(\s^Ny)/2\rfl}}{\widetilde{\Psi}_{N}^{(\s)}(y)+
  \frac{a^{\lfl\phi(\s^Ny)/2\rfl}-1}{a-1}}\ \ \text{by \ \ref{faAnchor}}\\ &>
 \frac{a-1}{1+\e(a-1)}.\ \ \CheckedBox\ \ge
 \end{align*}

\demo{ Proof that  $\overline{M}\le a-1$}
\

Fix $\e>0$.
\

For $n\ge 1\ \&\ y\in Y$, let  as in \ \ \ref{faBell}\ ,\  $n=\phi_{K_n}(y)+r_n(y)$, then
$$\Psi(\tau^n(y,0))=a^{R_n}\ \text{with}\ {R_n}=r_n(y)\wedge(\phi(\s^{K_n}y)-r_n(y))$$ whence

$$\Psi^{(\tau)}_{r_n}(\s^{K_n}y,0)=\sum_{k=0}^{r_n-1}a^{(k\wedge\phi(\s^{K_n}y)-k)}\ge\sum_{k=0}^{{R_n}-1}a^k=\tfrac{a^{R_n}-1}{a-1}.$$
Choose $n=n(y)\ge 1$ so large that
\begin{align*}\tag*{{\scriptsize\faBeer}}\label{faBeer}
 \tfrac{a-1}{\e\widetilde{\Psi}^{(\s)}_{K_n}(y)}<\tfrac{a-1}{1-\e}.
\end{align*}

Applying all this to\  \ \ref{faBank},
 \begin{align*}
M_n(\Psi,\tau)(y,0)&\le
 \frac{a^{R_n}}{ \widetilde{\Psi}_{K_n}^{(\s)}(y)+\frac{a^{R_n}-1}{a-1}}\\ &=
 \frac{a-1}{1-a^{-R_n}+a^{-R_n} \widetilde{\Psi}_{K_n}^{(\s)}(y)}\\ &
\le \tfrac{a-1}{1-\e}1_{[a^{-R_n}<\e]}+
\tfrac{a-1}{\e\widetilde{\Psi}_{K_n}^{(\s)}(y)}1_{[a^{-R_n}\ge\e]}\ \text{by\ \ref{faBeer}}\\ &
\lesssim\ \tfrac{a-1}{1-\e}.\ \ \ \CheckedBox
 \end{align*}
 \demo{Proof of Theorem 3.1}
 \

 Fix $f\in\mathcal{P}(\Bbb N)$ satisfying
$$\sum_{n\ge 1}nf(\{n\})<\infty\ \&\
\sum_{n\ge 1}a^nf(\{n\})=\infty\ \forall\ a>1.$$
 {\tt\large e.g.} any $f$ with $f(\{n\})\asymp \tfrac1{n^s}$ with $s>2$.

 \

Let $(\Om,m,S)$ be the the Markov shift with state space $\Bbb N$ and transition matrix $p:\Bbb N\x\Bbb N\to [0,1]$ given by
$$p_{s,t}=\begin{cases}& f_t\ \ \ s=1;\\ & 1\ \ \ s=t+1;\\ & 0\ \ \ \text{else.}
          \end{cases}$$
          As  on p.\pageref{renewal}, $(\Om,m,S)$ is isomorphic to the renewal process
          $(X,m,T)=(Y,p,\s)^\phi$  where
$$(Y,p,\s):=(\Bbb N^\Bbb Z,f^\Bbb Z,\text{\tt shift})$$ and $\phi:Y\to\Bbb N$ is defined
 by $\phi(y)=\phi((y_n:\ n\in\Bbb Z)):=y_0$, then $\Bbb E(\Phi)<\infty$.

\

Fix $t>0$ and define let $a=t+1$.

By construction, $(\sqrt{a}^{\phi\circ\s^n}:\ n\in\Bbb Z)$ are iid random variables with $\Bbb E(\sqrt{a}^\phi)=\infty$ and by Theorem 2.2,
$\mathfrak{e}(a^\Phi,\s)=\infty$ a.s..

Define $g=g_t:X\to\Bbb R_+$ by
$$g(y,t):=a^{t\wedge \phi(y)-t},$$ then, by Lemma 3.3
$$\mathfrak{e}(g,T)=t\ \text{a.s.}$$
To finish, we show that $g\circ\mathcal{c}:\Om\to\Bbb R_+$ is finitary (where $\mathcal{c}$ is as in \ref{faArrowsH} on p.\pageref{faArrowsH}).
\

Now
\begin{align*}
 g\circ\mathcal{c}(\om)&=g(\eta(\om),\mathcal{k}(\om))\\ &=
 a^{(-c_0(\om)+1)\wedge (\eta_0(\om)-(-c_0(\om)+1))}
\end{align*}
and $g\circ\mathcal{c}$ is finitary with $N(\om)=(-c_0(\om))\vee c_1(\om)$.\ \ \CheckedBox

\

For the proof of Theorem 3.2, we'll also need
\subsection*{  Dyadic ergodic stationary processes}
\

 Let $\Om:=\{0,1\}^\Bbb N,\ 
P:=\prod(\tfrac12,\tfrac12)$
\

The {\it dyadic odometer} $\tau:\Om:=\{0,1\}^\Bbb N\hookleftarrow$ is defined by

$$\tau(\om)=\tau(\undersetbrace{\ell-1\text{\rm-times}}\to{1,\dots,1},0,\om_{\ell+1},\dots):=
(\undersetbrace{\ell-1\text{\rm-times}}\to{0,\dots,0},1,\om_{\ell+1},\dots)$$
and $\tau(\mathbb{1}):=\mathbb{0}$.

Define $\ell:\Om\to\Bbb N\cup\{\infty\}$ by $\ell(\om):=\min\,\{n\ge 1:\ \om_n=0\}$.
\

Note that
\begin{align*}\tag*{\Bat}\label{Bat}\begin{split}
                                     & \{\tau^k\om|_{[1,n]}:\ 0\le k<2^n\}=\{0,1\}^n\ \forall\ n\ge 0;\\ &
                                     \forall\ n\ge 1,\ \exists\ 0\le k_n=k_n(\om)<2^{n-1}\ \text{s.t.}\ \ell(\tau^{k_n}\om)=n+\ell(\s^n\om)
                                    \end{split}\end{align*}

 where
$\s:\Om\to\Om$ is the shift: $\s(\om_1,\om_2,\dots):=(\om_2,\om_3,\dots)$.

\

A {\it dyadic stationary process} is a stationary process of form
$(\Om, P,\tau,\v)$ where

$\v(\om)=\b(\ell(\om)),\ \b:\Bbb
N\to\Bbb R_+,\ \b \uparrow$ and $\ell(\om):=\min\,\{n\ge 1:\ \om_n=0\}$.
\

For  a $\Bbb R_+$-ESP  $(\Om, P,\tau,\v)$ with
$\v=\b\circ\ell$,
\begin{align*}\tag*{\dsrailways}\label{dsrailways}
 \v_{2^n}(\om)&:=\sum_{k=0}^{2^n-1}\v(\tau^k\om)\ \ \overset{\text{\tiny\ref{Bat} }}=\
 \sum_{\e\in\{0,1\}^n\setminus\{\u 1\}}\b(\ell(\e))+\b(n+\ell\circ\s^n)\\ &=\sum_{k=1}^n2^{n-k}\b(k)+\b(n+\ell(\s^n\om)).
\end{align*}

\proclaim{Proposition 3.4}
\begin{align*}
 \tag*{\dsaeronautical}\label{dsaeronautical}\mathfrak{e}(a^\ell,\tau)=\infty\ \text{a.s.}\ \forall\ a\ge 2.
\end{align*}\endproclaim
\demo{Proof}
\ \ Fix $a\ge 2$ and write $\v:=a^\ell$ and $M_n:=\tfrac{\v\circ\tau^n}{\v_n}$.
\

For $n\ge 1,\ \om\in\Om$, let $k_n=k_n(\om)$ be as in \ref{Bat},\ then
$$\v_{k_n}\le \v_{2^n}-\v\circ\tau^{k_n}=2^n\sum_{k=1}^n(\tfrac{a}2)^k.$$
In case $a>2$,
$$\v_{k_n} \le Ca^n$$
for some fixed $C>0$ and $\forall\ n\ge 1$. Thus
\begin{align*}
M_{k_n}&=\tfrac{a^{n+\ell\circ\s^n}}{\v_{k_n}}\ge \tfrac{a^{n+\ell\circ\s^n}}{Ca^n}=\tfrac1C a^{\ell\circ\s^n}.
 \end{align*}
 To continue we claim that a.s.,
 \begin{align*}\tag*{\faBook}\label{faBook}\varlimsup_{n\to\infty}\ell\circ\s^n-\log_2n=\infty.
 \end{align*}
In particular,  $\varlimsup_{n\to\infty}\ell\circ\s^n=\infty$ a.s.
and  $\varlimsup_{n\to\infty}M_{k_n}=\infty$ a.s. when $a>2$.

In case $a=2$
$$\v_{k_n}\le\v_{2^n}-\v\circ\tau^{k_n}=n2^n$$
and
\begin{align*}
M_{k_n}&=\tfrac{2^{n+\ell\circ\s^n}}{\v_{k_n}}\ge \tfrac{2^{n+\ell\circ\s^n}}{n2^n}=2^{\ell\circ\s^n-\log_2n};\\ &
\overset{\text{\tiny\ref{faBook}}}\Rightarrow\ \varlimsup M_{k_n}=\infty\ \text{a.s.\ \CheckedBox}
 \end{align*}
 \demo{Proof of \ref{faBook}}

 We  show   $\ell\circ\s^n>\log_2n+r$ i.o. a.s. $\forall\ r\ge 1$.
 \

 To see this, fix $r\ge 1$, let $b_n\uparrow\infty$ be defined by $b_{n+1}=b_n+\kappa_n+r+1$ where
 $\kappa_n:=\lcl\log_2b_n\rcl$, then calculation shows that

\begin{align*}&\log_2 b_n\le \log_2n+\log_2\log_2n+o(1)\ \text{as}\ n\to\infty;
\end{align*}

 Now  let
 $A_n:=\{\om\in\Om:\ \om_k=1\ \forall\ b_n+1\le k\le b_n+\kappa_n+r\}$, then
\sbul $A_n\subset [\ell\circ\s^{b_n}\ge \log_2b_n+r]$;
\sbul $\{A_n:\ n\ge 1\}$ are independent (wrt $P$);
\sbul $P(A_n)=1/2^{r+\kappa_n}\gg \tfrac1{n\log_2n}$;
whence $\sum_{n\ge 1}P(A_n)=\infty$ and by the (classical) Borel-Cantelli lemma $\sum_{n\ge 1}1_{A_n}=\infty$ a.s.
Thus $\ell\circ\s^n>\log_2n+r$ i.o. a.s..\ \ \CheckedBox \ \ \ \ref{faBook} \ and hence \ref{dsaeronautical}\ \ when\ $a=2$.

\proclaim{Lemma 3.5}\ \ Fix $b\in\Bbb N_2$, let $(X_b,m_b,T_b):=(\Om,P,\tau)^{b\ell}$, and for $a>0$, let $\Psi_{a,b\ell}$ be as in the Main Lemma, then
$$\mathfrak{e}(\Psi_{a,b\ell},T_b)=a-1\ \text{a.s.}\ \forall\ a\ge 4^{\frac1b}.$$
\endproclaim\demo{Proof}\ \ Evidently $\Bbb E(b\ell)=2b$ and by Proposition 3.4,
$\mathfrak{e}(\sqrt{a}^{b\ell})=\infty$ a.s. $\forall\ a>4^{\frac1b}$. Thus, by the
Main Lemma,
$$\mathfrak{e}(\Psi_{r,b\ell},T_b)=a-1\ \text{a.s.}\ \ \CheckedBox$$

\

\

\proclaim{Lemma 3.6}\footnote{See also \cite[Corollary 5.6]{ORW}.}
\

\ For each $b\in\Bbb N\ \exists\ A_b\in\B(\Om)$ and an isomorphism
$\varpi_b:(A_b,P_{A_b},\tau_{A_b})\to (X_b,m_b,T_b)$.
\endproclaim\demo{Proof} \ \ For each $N\ge 1$,
$(\Om,P,\tau)^{2^N}\cong  (\Om,P,\tau)$.
\

Given $b\in\Bbb N$, choose $N\in\Bbb N$ so that $\Bbb E(b\ell)=2b<2^N$.
\

By \cite[Lemma 1.3]{ORW}, $\exists \ h:X_b\to\Bbb N$ measurable so that

$(\Om,P,\tau)\cong (X_b,m_b,T_b)^h$. The lemma follows from this.\ \CheckedBox

\demo{Proof of Theorem 3.2 {\rm (p.\pageref{3.2})}}
\

Fix $r>0$ and choose $b\in\Bbb N,\ b\ge 2$ so that $r\ge 4^{\frac1b}-1$.
\

Let
$(X_b,m_b,T_b)$ and $\Psi_{r,b\ell}:X_b\to\Bbb R_+$ be as in Lemma 3.5 so that $\mathfrak{e}(\Psi_{r,b\ell},T_b)=r$ a.s..

Let $(X,m,T)$ be an aperiodic, ergodic, probability preserving transformation.

By the Odometer Factor Proposition in \cite{AW2018}, $\exists$ $B\in\B(X)$ and a factor map $\pi:(B,m_B,T_B)\to (\Om,P,\tau)$ (the dyadic odometer).

Let $A_b\in\B(\Om)\ \&\ \varpi_b:A_b\to X_b$ be as in Lemma 3.6 and let
 $C:=\pi^{-1}A_b$. It follows that $\Pi:=\varpi_b\circ\pi|_{C}:(C,m_{C},T_{C})\to (X_b,m_b,T_b)$ is a factor map.

 Define $\psi:X\to [0,\infty)$ by
$\psi:=\Psi_{r,b\ell}\circ\Pi$ on $C$ and $\psi:=0$ on $X\setminus C$, then $\mathfrak{e}(\psi,T_C)=r$ a.s. on $C$ whence,
writing for a.e. $x\in C$,
$$M_n(T,x):=\tfrac{\psi(T^nx)}
{\sum_{k=0}^{n-1}\psi(T^kx)},$$\begin{align*}\mathfrak{e}(\psi,T)(x)&=
\varlimsup_{n\to\infty}M_n(T,x)\\ &=\varlimsup_{n\to\infty,\ T^nx\in C}M_n(T,x)\ \ \because M_n(T,x)=0\ \forall\ T^nx\notin C,\\ &=
 \varlimsup_{n\to\infty}M_n(T_C,x)\ \because\ \psi|_{X\setminus C}\equiv 0\\ &=\mathfrak{e}(\psi,T_C)(x)\ =\ r\ \ \text{a.s.}\ \CheckedBox
\end{align*}

\section*{\S4 Irrationality}
\subsection*{The Gauss map}
\

The Gauss map $G:\Bbb I\hookleftarrow$ is piecewise invertible with
inverse branches $\g_{[k]}:\Bbb I\to [k]:=[a=k]=(\tfrac1{k+1},\tfrac1k],\ \g_{[k]}(y)=\tfrac1{y+k}$.
\

Similarly, for  each $n\ge 1$, the inverse branches of $G^n:\Bbb I\hookleftarrow$ are
$\g_A:\Bbb I\to A$ where
$$A\in\a_n:=\{[a\circ G^k=a_k\ \forall\ 0\le k<n]:\ (a_0,a_1,\dots,a_{n-1})\in\Bbb N^n\}$$
 of form $\g_A:=\g_{[a_0]}\circ \g_{[a_1]}\circ\dots\circ \g_{[a_{n-1}]}$
 ($A=[a\circ G^k=a_k\ \forall\ 0\le k<n]$.
\

Writing, for $x\in \mathbb{I}\ \&\ n\in\Bbb N$, $x\in\a_n(x)\in\a_n$, we have
\begin{align*}
x & =\g_{\a_n(x)}(G^{n}x)
\\ &=\cfrac[r]{1|}{|a(x)}+\cfrac[r]{1|}{|a(Gx)}+\dots+\cfrac[r]{1|}{|a(G^{n-1}x) +G^{n}x}
\\ &\xrightarrow[n\to\infty]{}\cfrac[r]{1|}{|a_1}+\cfrac[r]{1|}{|a_2}+\dots+\cfrac[r]{1|}{|a_n} {\ \ +\ \dots}
\end{align*}
(where $a_n:=a(G^{n-1}x)$)  which latter is known as the {\it continued fraction expansion} of $x\in\mathbb{I}$.
\

The inverse to the continued fraction expansion is
$\frak b:\Bbb N^\Bbb N\to\mathbb{I}$  defined by
\begin{align*}\tag*{\dsarchitectural}\label{dsarchitectural}
\frak b(a_1,a_2,\dots) :=\cfrac[r]{1|}{|a_1}+\cfrac[r]{1|}{|a_2}+\dots+\cfrac[r]{1|}{|a_n} {\ \ +\dots}.
\end{align*} It is a homeomorphism $\mathfrak{b}:\Bbb N^\Bbb N\to\mathbb{I}$ conjugating
the Gauss map with the shift $S:\Bbb N^\Bbb N\hookleftarrow,\ \mathfrak{b}\circ S=G\circ\mathfrak{b}$.
\

Calculation shows that
$(\mathbb{I},m,G^2,\a_2)$ is an {\tt Adler map}, as in \cite{Fexp} satisfying
\begin{align*}&\tag{U}(G^2)^{\prime}\ge 4;\\ &
 \tag{A} \sup_{x\in \mathbb{I}}\tfrac{|(G^2)^{\prime\prime}(x)|}{(G^2)^{\prime}(x)^2}=2.
\end{align*}
It follows that
\begin{align*}
|\tfrac{\g_A^{\prime\prime}(x)}{\g_A^{\prime}(x)}|\le 4\ \forall\ n\ge 1,\ A\in\a_n,\ x\in \mathbb{I}.
\end{align*}
whence
\begin{align*} \tag{$\D$}\label{D}
|\g_A^{\prime}(x)|=e^{\pm 4}m(A)\ \forall\ n\ge 1,\ A\in\a_n,\ x\in \mathbb{I}.
\end{align*}
In particular, $m$ is a  Renyi measure for $G$.
\

Moreover by (\ref{D}),
$(\Bbb I,m,G,\{[a=m]:\ n\ge 1\})$ is a Gibbs-Markov map   and hence $d\mu(x)=\tfrac{dx}{\log 2(1+x)}$ is a {\tt c.f.} mixing measure for $G$ (see \cite{GM}).

\

\subsection*{ Convergents and denominators}\label{cgts}
\

The rest of this section is a collection of facts (from \cite{Kh} and \cite[\S4]{BillingsleyErg}) which we'll need in the sequel.
\

\par  Define the {\it convergents} $\tfrac{p_n}{q_n}$  ($p_n,\ q_n\in\Bbb Z_+,\ \gcd\,(p_n,q_n)=1$) of $x\in\mathbb{I}$  by
\begin{align*}
 \frac{p_n(x)}{q_n(x)}:=\cfrac[r]{1|}{|a(x)}+\cfrac[r]{1|}{|a(Gx)}+\dots+\cfrac[r]{1|}{|a(G^{n-1}x)}.
\end{align*}

\sbul The {\it principal denominators of $x$}  $q_n(x)$  are given by
$$q_{0}=1,\ q_{1}(x)=a(x),\ q_{n+1}(x)=a(G^nx)q_{n}(x)+q_{n-1}(x);$$
\sbul the numerators $p_n(x)$  are given by
$$p_{0}=0,\ p_{1}=1,\ p_{n+1}(x)=a(G^nx)p_{n}(x)+p_{n-1}(x).$$
It follows (inductively) that
\begin{align*}\tag*{\scriptsize\faTaxi}\label{faTaxi}
 q_n(x)\ge 2^{\frac{n-1}2}, \ p_n(x)=q_{n-1}(Gx)\ge 2^{\frac{n-2}2}\ \&\
 |x-\tfrac{p_n(x)}{q_n(x)}|<\tfrac1{q_n(x)q_{n+1}(x)}<\tfrac{\sqrt 2}{2^{n}}.
\end{align*}

Moreover:

\

\proclaim{4.1 Denominator lemma}\label{ 4.1}\ \ \cite[\S4]{BillingsleyErg},\ \cite{Kh}
\

\begin{align*}\tag*{\scriptsize\faRocket}\label{faRocket}|\log q_n(x)-\sum_{k=0}^{n-1}\log \tfrac1{G^k(x)}|\le \tfrac{2}{\sqrt{2}-1}\ \forall\ n\ge 1,\ x\in \mathbb{I}.
\end{align*}\endproclaim
\

It follows from  Birkhoff's theorem $\&$ \ref{faRocket} that if  $\mu\in\mathcal P(\mathbb{I})$ is $G$-invariant, ergodic, then
 \begin{align*}\tag*{\faPlane}\label{faPlane}\frac{\log q_n}n\ \underset{n\to\infty}\lra\ \int_\mathbb{I}\log\tfrac1xd\mu(x)\le\infty\ \ \ \mu-\text{a.s. }\ .\end{align*}\endproclaim
 \

 Also:
\proclaim{4.2 Proposition   \cite[\S4]{BillingsleyErg}, \cite[Th. 9 $\&$ 13]{Kh}}\label{4.2}
\

\begin{align*}\tag*{\dsmilitary}\label{dsmilitary}|x-\tfrac{p_n(x)}{q_n(x)}|=2^{\pm 1}\ \tfrac{G^n(x)}{q_n(x)^2}\ \forall\ n\ge 1,\ x\in \mathbb{I}.
\end{align*}\endproclaim
\proclaim{4.3 Corollary}\label{4.3}
\begin{align*}\tag*{\dsliterary}\label{dsliterary}
 m(\a_n(x))=(2M)^{\pm 1}\tfrac1{q_n(x)^2}\ \forall\ n\ge 1,\ x\in\Bbb I.
\end{align*}
\endproclaim\demo{Proof}

\begin{align*}
|x-\tfrac{p_n(x)}{q_n(x)}|&=|\g_{\a_n}(G^n(x))-\g_{\a_n(x)}(0)|\\ &=
G^n(x)|\g_{\a_n(x)}'(\th_n G_n(x))|\ \text{by Lagrange's theorem where}\ \th_n(x)\in [0,1]\\ &=
M^{\pm 1 }G^n(x)m(\a_n(x))\ \text{by \ref{Rp} on p\pageref{Rp}}
\end{align*}
and \ref{dsliterary} follows from \ref{dsmilitary} (p\pageref{dsmilitary}).\ \ \CheckedBox

\proclaim{4.4 Bugeaud's Lemma\ \ }\label{4.4}
\

\begin{align*}\tag*{\faKey}\label{faKey} \mathbbm{i}(x)=2+\mathbbm{e}((\log\tfrac1{G^nx}:\ n\ge 0))\ \forall\ x\in\mathbb{I}.
\end{align*}
\begin{align*}\tag*{\faPaw}\label{faPaw}\ \ \text{For $\mu\in\mathcal{P}(\mathbb{I})$ $G$-invariant,}\ \mathbbm{i}=2+\mathfrak{e}(\log a,G)\ \
\mu-\text{a.s..}
\end{align*}
\endproclaim
Note that \ref{faKey}  is a  version of
  \cite[Exercise E1]{Bugeaud}.
\demo{Proof of \ref{faKey}} \ \

Fix $x\in\Bbb I$. If  $x=\tfrac{\sqrt 5-1}2$, then
$a(G^nx)= 1\ \forall\ n\ge 0$ and $\mathbbm{e}((\log a(G^nx):\ n\ge 0))=\mathbbm{e}(\overline{0})=0$ and $\mathbbm{i}(x)=2$ (since $x$ is quadratic).
\

If $x\ne \tfrac{\sqrt 5-1}2$, $\exists\ \nu\ge 0,\ a(G^\nu x)\ge 2$ and for  $n\ge \nu,\ \sum_{k=0}^{n-1}\log a(G^kx)>0$.

Write $\widetilde{a}(x):=\tfrac1x$ and
$$M_n(x):=\frac{\log \widetilde{a}(G^nx)}{\sum_{k=0}^{n-1}\log \widetilde{a}(G^kx)},$$
then $\mathbbm{e}((\log \widetilde{a}(G^nx):\ n\ge 0))=\varlimsup_{n\to\infty}M_n(x)=:M(x)$.

\

We'll show that $M(x)=\mathbbm{i}(x)-2$ for $x\in\Bbb I$.
\

To this end, we show first that  $\sum_{n\ge 1}\log\widetilde{a}(G^n(x))=\infty$.

\

If $x\in\Bbb I,\ a(G^nx)\xrightarrow[n\to\infty]{}1$, then $\log\widetilde{a}(G^nx)\to\log\widetilde{a}(\tfrac{\sqrt 5-1}2)>0$
and $\sum_{n\ge 1}\log\widetilde{a}(G^n(x))=\infty$.
\

Otherwise, $\#\{n\ge 1:\ a(G^nx)\ge 2\}=\infty$ and
$$\sum_{n\ge 1}\log\widetilde{a}(G^n(x))\ge\log 2\#\{n\ge 1:\ a(G^nx)\ge 2\}=\infty.\ \CheckedBox$$
\

By \ref{dsmilitary} on p.\pageref{dsmilitary} , for $n\ge\nu\ \&\ \g>0$, we have
\begin{align*} q_n(x)^{2+\g}&|x-\tfrac{p_n(x)}{q_n(x)}| \asymp \frac{q_n(x)^{1+\g}}{q_{n+1}(x)}
\asymp \frac{q_n(x)^\g}{\widetilde{a}(G^nx)}
\\ &\asymp\exp[-(\log \widetilde{a}(G^nx)-\g\sum_{k=0}^{n-1}\log \widetilde{a}(G^kx))]\ \text{by \ref{faRocket} on  p.\pageref{faRocket}}\\ &=
\exp[(\sum_{k=0}^{n-1}\log \widetilde{a}(G^kx))(\g-M_n(x))]\\ &
\begin{cases}&\xrightarrow[n\to\infty]{}\infty\ \ \text{if}\ \ \ \ \g>M(x)
\\ &\to 0\ \text{along a subsequence if}\ \ \g<M(x).
\end{cases}
\end{align*}
Thus, $\mathbbm{i}(x)=M(x)+2$.\ \ \CheckedBox\ \ref{faKey}
\

\demo{Proof of \ref{faPaw}} By \ref{faKey}, $\mathbbm{i}=2+\mathfrak{e}(\log\widetilde{a},G)\ \mu$-a.s. and \ref{faPaw}
\ \ follows from Proposition 2.1 (below) since $|\log\widetilde{a}-\log a|\le 1$ on $\mathbb{I}$.\ \CheckedBox

\proclaim{4.5 Corollary}\label{4.5}
\

\sms{\rm (i)} If $\mu\in\mathcal{P}(\Bbb I)$ is so that $(\Bbb I,\mu,G,a)$ is {\tt c.f.} mixing,
then $\mu$-a.s. $x\in\Bbb I$ is Diophantine if $\Bbb E_\mu(\log a)<\infty$ and $\mu$-a.s. $x\in\Bbb I$ is Liouville if $\Bbb E_\mu(\log a)=\infty$;
\sms {\rm(ii)} For each $r\in \Bbb R_+,\ \exists\ p_r\in\mathcal{P}(\Om)$, $G$-invariant, ergodic
so that $\mathbbm{i}=2+r$  $p_r$-a.s..\endproclaim\demo{Proof}
Statement (i) {\sl[(ii)]}  follows from Proposition 2.1(b) and Theorem 2.2 {\sl[3.1]}.
\CheckedBox

\section*{\S5 Khinchin's dichotomy for weak Renyi processes of partial quotients}
\

\proclaim{  Borel-Cantelli Lemma for weak Renyi maps}\label{bcl}\ \  Suppose that $(\mathbb{I},m,T,\a)$ is a weak  Renyi map and let $n\ge 1,\ A_n\in\s(\a)$.
\par If $\sum_{k=1}^\infty m(A_k)=\infty$, and either {\rm (a)} $T$ is exact or {\rm (b)} $T$ is ergodic and $A_{n+1}\subset A_n$, then
$m(\varlimsup_{n\to\infty}T^{-n}A_n)=1$.
\endproclaim

\demo{Proof}
\par By \ref{wRp} (on p.\pageref{wRp}),   $\exists\ C>1$ such that
$$m(T^{-k}A_k\cap T^{-n}A_n)\le C m(T^{-k}A_k)m(T^{-n}A_n)\ \forall\ n\ne k.$$
Suppose that  $\sum_{k=1}^\infty m(A_k)=\infty$ and let
$$A_\infty:=[\sum_{k=1}^\infty 1_{A_k}\circ T^k=\infty]=\varlimsup_{n\to\infty}T^{-n}A_n.$$
\

By the  Erdos-Renyi Borel-Cantelli lemma (\cite{ErdosRenyi}\ $\&$/or  \cite[p.391]{Ren2})
$m(A_\infty)\ge\tfrac1C>0$.
\

In addition, $A_\infty\in\mathcal{T}(T)$ and $m(A_\infty)=1$ if $T$ is  exact.
\

Under assumption (b), $T^{-1}A_\infty\supseteq A_\infty$, whence $T^{-1}A_\infty=A_\infty\ \mod m$ and by ergodicity of $T$,
$m(A_\infty)=1$.
\ \ \CheckedBox

We'll call a measure $\mu\in\mathcal{P}(\mathbb{I})$  {\it doubling at $0$}\label{doubling} if
 \begin{align*}\tag*{\faCopy}\label{faCopy}
 \exists\ M>1,\ r_0>0\ \text{so that}\ \mu((0,2r))\le M\mu((0,r))\ \forall\ 0<r\le r_0.
 \end{align*}

\proclaim{5.1 Khinchin type dichotomy}\label{5.1}
\

Let $\mu\in\mathcal{P}(\mathbb{I})$ be an ergodic, weak Renyi measure for $G$ which is doubling at $0$.
 and let $f:\Bbb N\to\Bbb R_+$ be such that
$nf(n)\downarrow 0$ as $n\uparrow\ \infty$.

\sms{\rm (i)}\ \ If
$\sum_{n\ge 1}\tfrac{\mu((0,nf(n))}n<\infty$, then
$$\min_p\,|x-\tfrac{p}q|/\tfrac{f(q)}{q}\xrightarrow[q\to\infty]{}\infty\ \text{for\ \ $\mu$- a.e. $x\in\Bbb T$}.$$
\sms{\rm (ii)}\ \ If
$\Bbb E_\mu(\log a)<\infty$ and
$\sum_{n\ge 1}\tfrac{\mu((0,nf(n))}n=\infty$, then
$$\varliminf_{q\to\infty}\min_p\,|x-\tfrac{p}q|/\tfrac{f(q)}{q}=0\ \text{for $\mu$-a.e. $x\in\Bbb I$}.$$
\endproclaim

\proclaim{Lemma 5.2}\label{5.2}
\

Let $\mu\in\mathcal{P}(\mathbb{I})$ be an ergodic, weak Renyi measure for $G$ and let $f:\Bbb N\to\Bbb R_+$ be such that
$nf(n)\downarrow 0$ as $n\uparrow\ \infty$.
\sms{\rm (i)}\ \ If
$\sum_{n\ge 1}\tfrac{\mu((0,nf(n))}n<\infty$, then
for\ \ $\mu$- a.e. $x\in\Bbb I$,
$$\#\bigg\{\tfrac{p}q\in\Bbb Q:\ |x-\tfrac{p}q|<\tfrac{f(q)}{2q}\bigg\}<\infty.$$
\sms{\rm (ii)}\ \ If
$\Bbb E_\mu(\log a)<\infty$ and
$\sum_{n\ge 1}\tfrac{\mu((0,nf(n))}n=\infty$, then
for\ \ $\mu$- a.e. $x\in\Bbb I$,
$$\#\bigg\{\tfrac{p}q\in\Bbb Q:\ |x-\tfrac{p}q|<\tfrac{f(q)}q\bigg\}=\infty.$$
\endproclaim

\proclaim{5.3 Remark}\label{5.3}\ \ Let $f:\Bbb R_+\to\Bbb R_+$ be such that
$xf(x)\downarrow 0$ as $x\uparrow\ \infty$.
\

Define $h:[1,\infty)\to[\tfrac1{f(1)},\infty)$ by
$h(x):=\tfrac1{xf(x)}$ and let $g=h^{-1}:[\tfrac1{f(1)},\infty)\to[1,\infty)$ then
$\Bbb E_\mu(\log g\circ a)<\infty$
iff $\sum_{n\ge 1}\tfrac{\mu((0,nf(n)))}n<\infty$.\endproclaim
\

\demo{Proof of Remark 5.3}
Fix $\kappa>1$, then $\Bbb E_\mu(\log g\circ a)<\infty$ iff
\begin{align*}\infty>&\sum_{n\ge 1}\mu([\log g\circ a>n\log\kappa])=
\sum_{n\ge 1}\mu([g\circ a>\kappa^n])\\ &\asymp
\sum_{n\ge 1}\tfrac{\mu([g\circ a>n])}n\ \text{by condensation,}\\ &=
\sum_{n\ge 1}\tfrac{\mu([a>g^{-1}(n)])}n=
\sum_{n\ge 1}\tfrac{\mu((0,\tfrac1{g^{-1}(n)})}n\\ &=
\sum_{n\ge 1}\tfrac{\mu((0,nf(n))}n.\ \ \CheckedBox
\end{align*}

In particular, with $f(x)=\tfrac1{x^{1+s}}$  ($s>0$), we have $g(x)=x^{\frac1s}$
and $\log g\circ a=\tfrac1s\log a$, whence
\begin{align*}\tag*{\faLifeRing}\label{faLifeRing}\Bbb E_\mu(\log a)<\infty\ \iff\ \sum_{n\ge 1}\tfrac{\mu((0,\frac1{n^s})}n<\infty\ \ \text{ for some (hence all) $s>0$}.\end{align*}
\demo{Proof of Lemma 5.2(i)}\ \
\

 By \ref{dsmilitary} on p.\pageref{dsmilitary},  we have  that
\begin{align*}|x-\tfrac{p_n(x)}{q_n(x)}|\ge\tfrac{G^n(x)}{2q_n(x)^2}\ \forall\ n\ge 1,\ x\in\mathbb{I}.
\end{align*}

Fix $1<\kappa<\exp[\int_\Om\log\tfrac1x d\mu(x)]$.   By condensation,

$\sum_{n\ge 1}\mu([0,\kappa^nf(\kappa^n)])<\infty$  and for
$\mu$-a.e. $x\in \mathbb{I},\ \exists\ N(x)$ so that
$$G^n(x)\ge\kappa^nf(\kappa^n)\ \forall\ n\ge N(x).$$
Moreover, by   \ref{faPlane} on p.\pageref{faPlane}, we can ensure that for $\mu$-a.s. $x\in\mathbb{I},\ \exists\
N_1(x)>N(x)$ so that in addition, $\forall\ n>N_1(x)$:
$$q_n(x)\ge \kappa^n\ \&\ \text{hence also}\ \kappa^nf(\kappa^n)\ge q_n(x)f(q_n(x)).$$
Thus, for $\mu$-a.s. $x\in \mathbb{I},\ n\ge N_1(x)$,

\begin{align*}\tag*{\dschemical}\label{dschemical}|x-\tfrac{p_n(x)}{q_n(x)}|\ge\tfrac{G^n(x)}{2q_n(x)^2}\ge
\tfrac{\kappa^nf(\kappa^n)}{2q_n(x)^2}\ge \tfrac{q_n(x)f(q_n(x))}{2q_n(x)^2}=\tfrac{f(q_n(x))}{2q_n(x)}.
\end{align*}
Lastly, if $|x-\tfrac{p}q|<\tfrac{f(q)}{2q}$ and $q$ is large enough so that
$\tfrac{f(q)}{2q}<\tfrac1{2q^2}$, then by Legendre's theorem
(see e.g. \cite[Theorem 5C]{WSchmidt}), $q=q_n(x)$ (some $n\ge 1$) and \ref{dschemical} applies contradicting $|x-\tfrac{p}q|<\tfrac{f(q)}{2q}$.
\ \ \CheckedBox\ {\rm (i)
\demo{Proof of Lemma 5.2(ii)}
 \

 We'll prove under the assumptions, that
 for $\mu$-a.s. $x\in \mathbb{I}$,
\begin{align*}
\#\bigg\{n\in\Bbb N:\ |x-\tfrac{p_n(x)}{q_n(x)}|<\tfrac{f(q_n(x))}{q_n(x)}\bigg\}=\infty.\end{align*}

To this end, fix $\kappa>\exp[\int_\Bbb I\log\tfrac1xd\mu(x)]$.
\

By condensation,
$\sum_{n\ge 1}\mu([a>\tfrac1{\kappa^nf(\kappa^n)}])=\infty$ and by the Borel-Cantelli lemma under assumption (b)   (on p.\pageref{bcl})
for $\mu$- a.s. $ x\in \mathbb{I}$,
$$\mu(\{x\in \mathbb{I}:\ \#\{n\ge 1:\ G^{n}x<\kappa^nf(\kappa^n)\}=\infty\}).$$
 \

 By  \ref{faPlane} on p.\pageref{faPlane},  for $\mu$-a.e. $x\in \mathbb{I}$ , \ $\#\{n\ge 1:\   q_n( x)\ge\kappa^n\}<\infty$ whence
\  $\# K(x)=\infty$ where
$$K(x):=\{n\ge 1:\ q_n(x)<\kappa^n\ \&\ G^{n}x<\kappa^nf(\kappa^n)\}.$$
For $n\in K(x)$, we have
\begin{align*} |x-\frac{p_n(x)}{q_n(x)}| &<\frac1{q_n(x)q_{n+1}(x)}<
\frac1{a(G^nx)q_n(x)^2}<\frac{\kappa^nf(\kappa^n)}{q_n(x)^2}\\ &
\le\frac{q_n(x)f(q_n(x))}{q_n(x)^2}\ \because\ kf(k)\downarrow\ \&\
q_n(x)<\kappa^n\\ &=
\frac{f(q_n(x))}{q_n(x)}.\ \ \CheckedBox\ {\rm (ii)}\end{align*}
\demo{Proof of Theorem 5.1}\ \ By the doubling property,
$$\sum_{n\ge 1}\tfrac{\mu((0,nf(n))}n\ \ \underset{\text{\large $=$}}{^{\text{\large $<$}}}\ \ \infty\ \iff\ \ \sum_{n\ge 1}\tfrac{\mu((0,cnf(n))}n\ \ \underset{\text{\large $=$}}{^{\text{\large $<$}}}\ \ \infty\ \forall\ c>0$$
so Lemma 5.2 holds for each $f_c:=cf$ ($c>0$).
\

Theorem 5.1 follows from this.\ \CheckedBox

\subsection*{Ahlfors-regular, Gauss-invariant  measures}
\

Consider the full shift $(X_K:=K^\Bbb N,S)$ where $K\subset\Bbb N$ is infinite and $S:K^\Bbb N\hookleftarrow$ is the shift. Let $Y_K:=\mathfrak{b}(X_K)\subset\mathbb{I}$ where $\frak b:\Bbb N^\Bbb N\to\Bbb I$ is as in
\ref{dsarchitectural} on p. \pageref{dsarchitectural}.
\

By \cite[Theorem 7.1]{FSU}, for each $h\in (0,1],\ \exists\ K=K(h)\subset\Bbb N$ infinite
so that the Hausdorff dimension of $Y_K$ is $h$; and so that $\mu_{K}\in\mathcal{P}(Y_{K})$, the restriction of the Hausdorff measure with gauge function $t\mapsto t^h$ to $Y_{K}$ is $h$-{\it Ahlfors-regular} in the sense that $\exists\ c>1$ so that
\begin{align*}\tag*{\dsmedical}\label{dsmedical}
\mu_K((x-\e,x+\e))=c^{\pm 1}\e^h\ \forall\ x\in\text{\tt Spt}\,\mu_K,\ \e>0\ \text{small}.
\end{align*}
\proclaim{5.4 Corollary\ \ (\cite[Theorem 6.1]{FSU})}\label{5.4}
 \

 Let $h\in (0,1]\ \&\ K\subset\Bbb N$ be infinite and let $\mu_K\in\mathcal{P}(Y_K)$ satisfy \ \ref{dsmedical} with parameter $h$, then
 $\Bbb E_{\mu_K}(\log a)<\infty$ and for $f:\Bbb N\to\Bbb R_+,\ nf(n)\downarrow$,
\begin{align*}\tag*{\faLock}\label{faLock}\min\,\{|x-\tfrac{p}q|:\ p\in\Bbb N\}\underset{q\to\infty}\gg \tfrac{f(q)}q\  \text{\small\rm for\ $\mu_K$-a.s.}\ x\in\mathbb{I}\ \text{\small\rm iff}\
\sum_{n\ge 1}\tfrac{f(n)^h}{n^{1-h}}<\infty.
\end{align*}\endproclaim
\f{\bf Remark}
\

As shown in \cite{BHZ}, in contrast to this, {\tt self similar measures} (which are also Ahlfors regular) satisfy \ref{faLock} with $h=1$, whatever their dimension $h\in (0,1]$.

\demo{Proof}
Since
$$GY_K=G\circ\mathfrak{b}(X_K)=\mathfrak{b}\circ S(X_K)=\mathfrak{b}(X_K)=Y_K,$$
it follows from \ref{dsmedical} (p.\pageref{dsmedical}) via Besicovitch's differentiation theorem (see e.g. \cite[Chapter 2]{Mattila})\ that for $n\ge 1,\ \mu_K\circ G^n\ll\mu_K$ with
\begin{align*}
 \tag*{\dsagricultural}\label{dsagricultural}\frac{d\mu_K\circ G^n}{d\mu_K}=
 c_K^{\pm 1}(|G^{n\prime}|)^h\ \mu_K-\text{a.s..}
\end{align*}
For $n\ge 1$, let
$$\b_n:=\{A\in\a_n:\ \mu_K(A)>0\},$$ then for $A\in\b_n$, $\mu_K$-a.s.,
\begin{align*}
 \tfrac{d\mu_K\circ\g_A}{d\mu_K}&=(\tfrac{d\mu_K\circ G^n}{d\mu_K}\circ\g_A)^{-1}\\ &=
 c^{\pm 1}|G^{n\prime}\circ\g_A|^{-h}\\ &=c^{\pm 1}|\g_A'|^h\\ &=M^{\pm 1}m(A)^h
 \ \text{by \ref{D} on p.\pageref{D}}
\end{align*}
where $M=ce^{4h}$.

Moreover
$$\mu_K(A)=\int_{\Bbb I}\tfrac{d\mu_K\circ\g_A}{d\mu_K}d\mu_K=M^{\pm 1}m(A)^h$$
with the conclusion that
$$\tfrac{d\mu_K\circ\g_A}{d\mu_K}=M^{\pm 2}\mu_K(A).$$
\

  By \cite{RenRep} $\exists\ P_K\in\mathcal{P}(Y_K)$,\ $P_K\sim\mu_K$ so that $P_K\circ G^{-1}=P_K$ and so that
$\log\tfrac{dP_K}{d\mu_K}\in L^\infty(\mu_K)$.
\

Thus  $(Y_K,P_K,G,\a\})$ has the  Renyi property.
\

Since $K$ is infinite,\ $0\in\text{\tt Spt}\,\mu_K$ and by \ref{dsmedical} (p.\pageref{dsmedical}),  $\mu_K((0,y))=c_K^{\pm 1}y^h\ \forall\ y>0$ small and in particular, $\mu_K$ is doubling at $0$.

By \ref{faLifeRing} on p.\pageref{faLifeRing},\ $\Bbb E_{\mu_K}(\log a)<\infty$.

Thus, \ref{faLock} follows from Theorem 5.1 (p.\pageref{5.1}).\ \CheckedBox\
\

\end{document}